\documentclass[12pt,reqno]{article}

\usepackage[usenames]{color}
\usepackage{amssymb}
\usepackage{graphicx}
\usepackage{amscd}

\usepackage[colorlinks=true,
linkcolor=webgreen,
filecolor=webbrown,
citecolor=webgreen]{hyperref}

\definecolor{webgreen}{rgb}{0,.5,0}
\definecolor{webbrown}{rgb}{.6,0,0}

\usepackage{color}
\usepackage{fullpage}
\usepackage{float}

\usepackage{graphics,amsmath,amssymb}
\usepackage{amsthm}
\usepackage{amsfonts}
\usepackage{latexsym}
\usepackage{epsf}

\begin{document}

\begin{center}
\epsfxsize=4in
%\leavevmode\epsffile{logo129.eps}
\end{center}

\theoremstyle{plain}
\newtheorem{theorem}{Theorem}
\newtheorem{corollary}[theorem]{Corollary}
\newtheorem{lemma}[theorem]{Lemma}
\newtheorem{proposition}[theorem]{Proposition}

\theoremstyle{definition}
\newtheorem{definition}[theorem]{Definition}
\newtheorem{example}[theorem]{Example}
\newtheorem{conjecture}[theorem]{Conjecture}

\theoremstyle{remark}
\newtheorem{remark}[theorem]{Remark}

\begin{center}
\vskip 1cm{\LARGE\bf Formulas for Odd Zeta Values \vskip .1in 
and Powers of $\pi$ }
\vskip 1cm
\large Marc Chamberland
and Patrick Lopatto\\
Department of Mathematics and Statistics\\
Grinnell College\\
Grinnell, IA 50112\\
USA \\
\href{mailto:chamberl@math.grinnell.edu}{\tt chamberl@math.grinnell.edu} \\
\end{center}

\vskip .2 in

\begin{abstract}
Plouffe conjectured rapidly converging series
formulas for $\pi^{2n+1}$ and $\zeta (2n+1)$
for small values of $n$. We find the general pattern
for all nonnegative integer values of $n$ and offer a proof.
\end{abstract}

\makeatletter
\newcommand{\rmnum}[1]{\romannumeral #1}
\newcommand{\Rmnum}[1]{\expandafter\@slowromancap\romannumeral #1@}

\newcommand{\qdp}{\hspace*{1.5mm}}%
\newcommand{\qqdp}{\hspace*{2.5mm}}
\newcommand{\xqdp}{\hspace*{5.0mm}}
\newcommand{\xxqdp}{\hspace*{10mm}}
\newcommand{\xquad}{\quad\qquad\quad}
\newcommand{\qdn}{\hspace*{-1.5mm}}
\newcommand{\qqdn}{\hspace*{-2.5mm}}
\newcommand{\xqdn}{\hspace*{-5.0mm}}
\newcommand{\xxqdn}{\hspace*{-10mm}}

%%%%%%%%%%%%%%%%%%%%%%%%%%%%%%%%%%%%%

\newcommand{\fns}{\footnotesize}
\newcommand{\dst}{\displaystyle}
\newcommand{\sst}{\scriptstyle}
\newcommand{\sss}{\scriptscriptstyle}

%%%%%%%%%%%%%%%%%%%%%%%%%%%%%%%%%%%%%

\newcommand{\+}{&\qqdn}%
\newcommand{\zero}{&\qdn&\xqdn}
\newcommand{\tagg}[2]{\addtocounter{equation}{#1}
        \tag{\theequation#2}}%%%\tag{}\tag*{}

%%%%%%%%%%%%%%%%%%%%%%%%%%%%%%%%%%%%%%%%%%%%%%%%%%%%%%
%DDDDDDDDDDDDDDDDDDDDDDDDDDDDDDDDDDDDDDDDDDDDDDDDDDDD%
%%%%%%%%%%%%%%%% symbols and noations %%%%%%%%%%%%%%%%

\newcommand{\la}{\leftarrow}%
\newcommand{\ra}{\rightarrow}
\newcommand{\La}{\Leftarrow}%
\newcommand{\Ra}{\Rightarrow}
\newcommand{\lla}{\longleftarrow}%
\newcommand{\lra}{\longrightarrow}
\newcommand{\Lla}{\Longleftarrow}%
\newcommand{\Lra}{\Longrightarrow}
\newcommand{\mto}{\mapsto}
\newcommand{\rmap}{\longmapsto}%%%
\newcommand{\ltor}{\rightleftharpoons}

\newcommand{\diam}{\diamond}
\newcommand{\ang}[1]{\langle{#1}\rangle}
\newcommand{\res}[2]{\mbox{Res}\left[#1\right]_{#2}}

\newcommand{\mb}[1]{\mathbb{#1}}
\newcommand{\mc}[1]{\mathcal{#1}}
\newcommand{\sbs}[1]{_{\substack{#1}}}
\newcommand{\sps}[1]{^{\substack{#1}}}

\newcommand{\dvd}[3]
           {\sum_{i=1}^n {#1}
            \prod\sbs{j=1\\j\not=i}^n
            \frac{x_i^{\!#2}}{x_i-x_j}{#3}}
\newcommand{\mdvd}[2]
           {\sum\sbs{\Lambda\subset[n]\\
                    |\Lambda|=p}{#1}
            \prod\sbs{i\in\Lambda\\j\not\in\Lambda}
            \frac{x_i}{x_i-x_j}{#2}}

\newcommand{\fat}[3]{\left(#1;#2\right)_{#3}}
\newcommand{\ffat}[3]{\left[#1;#2\right]_{#3}}
\newcommand{\fnk}[3]{\left[\qdn\ba{#1}#2\\#3\ea\qdn\right]}
\newcommand{\ffnk}[4]{\left[\qdn\ba{#1}#3\\#4\ea{\!\Big|\:#2}\right]}
\newcommand{\fnkl}[4]{\left[{\qdn\ba{#1}#3\\#4\ea{#2}}\right.}
\newcommand{\fnkr}[4]{\left.{\ba{#1}#3\\[3.5mm]#4\ea}{\Big|\:#2}\right.}

\newcommand{\binm}{\binom}
\newcommand{\sbnm}[2]{\Bigl(\!\ba{c}\!#1\!\\#2\ea\!\Bigr)}
\newcommand{\tbnm}[2]{\tbinom{#1}{#2}}
\newcommand{\zbnm}[2]{\genfrac{(}{)}{0mm}{2}{#1}{#2}}
\newcommand{\binq}[2]{\genfrac{[}{]}{0mm}{0}{#1}{#2}}
\newcommand{\sbnq}[2]{\Bigl[\!\ba{c}\!#1\!\\#2\ea\!\Bigr]}
\newcommand{\tbnq}[2]{\genfrac{[}{]}{0mm}{1}{#1}{#2}}
\newcommand{\zfrac}[2]{\genfrac{}{}{}{2}{#1}{#2}}
\newcommand{\double}[2]{\genfrac{}{}{0mm}{1}{#1}{#2}}
\newcommand{\cinque}{\mult{c}{:\\[-2.5mm]\vdots}}

\newcommand{\nnm}{\nonumber}
\newcommand{\be}{\begin{equation}}
\newcommand{\ee}{\end{equation}}
\newcommand{\ba}{\begin{array}}
\newcommand{\ea}{\end{array}}
\newcommand{\bmn}{\begin{eqnarray}}
\newcommand{\emn}{\end{eqnarray}}
\newcommand{\bnm}{\begin{eqnarray*}}
\newcommand{\enm}{\end{eqnarray*}}
\newcommand{\bln}{\begin{subequations}}
\newcommand{\eln}{\end{subequations}}
\newcommand{\blm}{\bln\bmn}
\newcommand{\elm}{\emn\eln}

\newcommand{\eqn}[1]{\begin{equation}#1\end{equation}}
\newcommand{\xalignx}[1]{\begin{align}#1%%%%%%%%%%%%%%
            \end{align}}     %%%align=auto-space%%%%%%
\newcommand{\xalignz}[1]{\begin{align*}#1%%%%%%%%%%%%%
            \end{align*}}    %%%align=auto-space%%%%%%
\newcommand{\xflalignx}[1]{\begin{flalign}#1%%%%%%%%%%
            \end{flalign}}   %%%flalign=max-space%%%%%
\newcommand{\xflalignz}[1]{\begin{flalign*}#1%%%%%%%%%
            \end{flalign*}}  %%flalign=max-space%%%%%%
\newcommand{\xalignatx}[2]{\begin{alignat}{#1}#2%%%%%%
            \end{alignat}}   %%%alignat=no-space%%%%%%
\newcommand{\xalignatz}[2]{\begin{alignat*}{#1}#2%%%%%
            \end{alignat*}}  %%%alignat=no-space%%%%%%

\newcommand{\xsplitx}[1]{\begin{split}#1\end{split}}%%
\newcommand{\xalignedx}[1]{\begin{aligned}#1%%%%%%%%%%
            \end{aligned}}  %%%aligned in equation%%%%
\newcommand{\xalignedatx}[2]{\begin{alignedat}{#1}#2%%
            \end{alignedat}}%%%alignedat in equation%%
\newcommand{\xgatherx}[1]{\begin{gather}#1%%%%%%%%%%%%
            \end{gather}}   %%%%%gather senza  & %%%%%
\newcommand{\xgatherz}[1]{\begin{gather*}#1%%%%%%%%%%%
           \end{gather*}}   %%%%%gather senza & %%%%%%
\newcommand{\xgatheredx}[1]{\begin{gathered}#1%%%%%%%%
            \end{gathered}}%gathered=no& in equation%%
\newcommand{\mult}[2]{\begin{array}{#1}#2\end{array}}%
\newcommand{\xmultx}[1]{\begin{multline}#1%%%%%%%%%%%%
            \end{multline}} %%%multline senza & %%%%%%
\newcommand{\xmultz}[1]{\begin{multline*}#1%%%%%%%%%%%
            \end{multline*}}%%%multline senza & %%%%%%

%%%%%%%%%%%%%%%%%%%%%%%%%%%%%%%%%%%%%%%%%%%%%%%%%%%%%%
%FFFFFFFFFFFFFFFFFFFFFFFFFFFFFFFFFFFFFFFFFFFFFFFFFFFF%
%%%%%%%%%%%%%%%%%%%%% positions %%%%%%%%%%%%%%%%%%%%%%

\newcommand{\centro}[1]
           {\begin{center}#1\end{center}}
\newcommand{\sinistra}[1]
           {\begin{flushleft}#1\end{flushleft}}
\newcommand{\destra}[1]
           {\begin{flushright}#1\end{flushright}}
\newcommand{\centerbox}[2]{\centro{\begin{tabular}{|c|}\hline
\parbox{#1}{\mbox{}\\[2.5mm]#2\mbox{}\\}\\\hline\end{tabular}}}
\newcommand{\centrobmp}[3]{\vspace*{1mm}
        \centerbmp{#1}{#2}{#3}\vspace*{1mm}}

\newcommand{\sav}[1]{\newsavebox{#1}%
 \sbox{#1}{\theequation}\\[-2mm]}%=\newline
\newcommand{\usa}[1]{(\usebox{#1})}%{#1}={\some-nome}

\newcommand{\summary}[2]{\begin{center}\parbox{#1}
{{\sc\bf Summary}:\,{\small\it#2}}\end{center}}
\newcommand{\riassunto}[2]{\begin{center}\parbox{#1}
{{\sc\bf Riassunto}:\,{\small\it#2}}\end{center}}
\newcommand{\grazia}[2]{\begin{center}\parbox{#1}
{{\sc\bf Ringraziamento}:\,{\small\it#2}}\end{center}}
\newcommand{\thank}[2]{\begin{center}\parbox{#1}
{{\sc\bf Acknowledgement}:\,{\small\it#2}}\end{center}}

\newcommand{\bbtm}[5]{\bibitem{kn:#1}{#2,}~{#3,}~\emph{#4}~{#5.}}
\newcommand{\bbtmn}[4]{\bibitem{kn:#1}{#2,}~\emph{#3,}~{#4.}}
\newcommand{\bbtmm}[4]{\bibitem{kn:#1}{#2,}~{#4.}}
\newcommand{\cito}[1]{\cite{kn:#1}}
\newcommand{\citu}[2]{\cite[#2]{kn:#1}}

\newcommand{\mr}{\textbf{MR}\:}
\newcommand{\zbl}{\textbf{Zbl}\:}

\newcommand{\nota}[2]{\centerbox{#1}{\textbf{Achtung}\quad#2}}

\newcommand{\graph}[3]{\begin{picture}(0,0)
         \put(#1,#2){#3}\end{picture}}
\newcommand{\fwd}{\mbox{$\bigtriangleup\qdn\!\cdot\,\:$}}
\newcommand{\bwd}{\bigtriangledown}

\section{Introduction} \label{sec-intro}
\setcounter{equation}{0}

It took nearly one hundred years for the Basel Problem --- finding a
closed form solution to $\sum_{k=1}^\infty 1/k^2$ --- to see a solution.
Euler solved this in 1735 and essentially solved the problem where
the power of two is replaced with any even power. This formula
is now usually written as
$$
\zeta (2n ) = (-1)^{n+1} \frac{ B_{2n}(2\pi)^{2n} }{ 2(2n)! } ,
$$
where $\zeta(s)$ is the Riemann zeta function and $B_k$ is
the $k$-th Bernoulli number, defined by the
generating function
$$
\frac{x}{e^x - 1} = \sum_{n=0}^\infty \frac{B_n x^n}{n!}, ~~~|x|<2\pi,
$$
whose first few values are $0,-1/2, 1/6, 0, -1/30, \dots$.
However, finding a closed form for
$\zeta(2n+1)$ has remained an open problem.
Only in 1979 did Ap\'ery show that $\zeta (3)$
is irrational. His proof involved the snappy acceleration
$$
\zeta (3) = \frac{5}{2}\sum_{n=1}^\infty
\frac{ (-1)^{n-1} }{ n^3 {2n\choose n}} .
$$
This tidy formula does not generalize to the other 
odd zeta values, but other representations, such as nested
sums or integrals, have been well studied.
The hunt for a clean result like Euler's
has largely been abandoned, leaving researchers with the goal
of finding formulas which either converge quickly or have an elegant form.

Following his success in discovering a new formula for $\pi$, 
Simon Plouffe \cite{Plouffe} conjectured several identities
which relate either $\pi^m$ or $\zeta (m)$ to three
infinite series. Letting
$$
S_n (r) = \sum_{k=1}^\infty \frac{1}{k^n (e^{\pi rk} - 1)},
$$
the first few examples are\footnote{There is a typographical error in the sign of the last coefficient in the formula for $\pi^5$ in \cite{Plouffe}, which is corrected here.}
\begin{eqnarray*}
\pi & = & 72 S_1 (1) - 96 S_1 (2) + 24 S_1 (4) \\
\pi^3 & = & 720 S_3 (1) - 900 S_3 (2) + 180 S_3 (4) \\
\pi^5 & = & 7056 S_5 (1) - 6993 S_5 (2)  -  63 S_5 (4) \\
\pi^7 & = & \frac{907200}{13} S_7 (1)-70875 S_7 (2)+\frac{14175}{13} S_7 (4).
\end{eqnarray*}
and
\begin{eqnarray*}
\zeta (3) & = & 28 S_3 (1) - 37 S_3 (2) + 7 S_3 (4) \\
\zeta (5) & = & 24 S_5 (1) - \frac{259}{10} S_5 (2)  -  \frac{1}{10} S_5 (4) \\
\zeta (7) & = & \frac{304}{13} S_7 (1) - \frac{103}{4} S_7 (2) 
+ \frac{19}{52} S_7 (4)  .
\end{eqnarray*}
Plouffe conjectured these formulas by first assuming, for example,
that there exist constants $a$, $b$, and $c$ such that
$$
\pi = a S_1 (1) + b S_1 (2) + c S_1 (4).
$$
By obtaining accurate approximations of each the three series, 
he wrote some computer code to postulate rational values for
$a,b,c$. Today, such integer relation algorithms have been 
used to discover many formulas. The widely used
PSLQ algorithm, developed by Ferguson and Bailey \cite{FB}, 
is implemented in Maple.
The following Maple code solves the above problem:
\begin{verbatim}
> with(IntegerRelations):
> Digits := 100;
> S := r -> sum( 1/k/( exp(Pi*r*k)-1 ), k=1..infinity );
> PSLQ( [ Pi, S(1), S(2), S(4) ] );
\end{verbatim}
The PSLQ command returns the vector $[-1,72,-96, 24]$, producing
the first formula.

While the computer can be used to conjecture the coefficients
for a specific power, finding the general sequences of rationals
has remained an open problem. This note finds these sequences and
offers formal proofs.

\section{ Exact Formulas }

While it does not seem that $\zeta(2n+1)$ is a rational multiple
of $\pi^{2n+1}$, a result in Ramanujan's notebooks gives a 
relationship with rapidly convergent infinite series. See Entry 21(i) in Chapter 14 of \cite{Berndt}, and see \cite{berndt2017ramanujan} for additional commentary.
% We define
% \begin{equation*}
% S_{n}(r)=\sum_{k=1}^{\infty}\frac{1}{k^{n}(e^{r k}-1)}
% \end{equation*}
% and let $B_j$ denote the $j^{th}$ Bernoulli number.

\begin{theorem}[Ramanujan]
\label{Rama_thm}
If $\alpha>0$, $\beta>0$, and $\alpha \beta = \pi^2$, then
% \begin{multline}
% \alpha^{-n}\left\{\frac{1}{2}\zeta(2n+1)+S_{2n+1}(2\alpha) \right\} = \\(-\beta)^{-n}\left\{\frac{1}{2}\zeta(2n+1)+S_{2n+1}(2\beta) \right\} - 4^n\sum_{k=0}^{n+1}(-1)^k\frac{B_{2k} B_{2n+2-2k}}{(2k)!(2n+2-2k)!}\alpha^{n+1-k}\beta^k
% \end{multline}

\begin{eqnarray*}
\lefteqn{
\alpha^{-n}\left\{\frac{1}{2}\zeta(2n+1)+S_{2n+1}(2\alpha/\pi) \right\} =  }\\
& & 
(-\beta)^{-n}\left\{\frac{1}{2}\zeta(2n+1)+S_{2n+1}(2\beta/\pi) \right\} 
- 4^n\sum_{k=0}^{n+1}(-1)^k\frac{B_{2k} B_{2n+2-2k}}{(2k)!(2n+2-2k)!}\alpha^{n+1-k}\beta^k .
\end{eqnarray*}

\end{theorem}
Using $\alpha=\beta=\pi$ in Proposition \ref{Rama_thm} and defining
$$
F_{n}=\sum_{k=0}^{n+1}(-1)^k\frac{B_{2k} B_{2n+2-2k}}{(2k)!(2n+2-2k)!} ,
$$
we have
$$
\left( \pi^{-n} - (-\pi)^{-n} \right)
\left(\frac{1}{2}\zeta(2n+1)+S_{2n+1}(2) \right) = 
- 4^n\pi^{n+1}F_n .
$$

To find formulas for the odd zeta values and powers of $\pi$, 
we will divide these into two classes: $\zeta (4m-1)$ and
$\zeta (4m+1)$. Such distinctions can be seen in other 
studies; see \cite[pp. 137--139]{BB}.

First, we find the formulas for $\pi^{4m-1}$ and $\zeta(4m-1)$. 
If $n$ is odd, then

\begin{equation}
\label{eqn_(3)}
\frac{1}{2}\zeta(2n+1)+S_{2n+1}(2)  = \frac{- 4^n}{2}\pi^{2n+1}F_n.
\end{equation}
Using $\alpha=\pi/2$ and $\beta=2\pi$ in 
Theorem \ref{Rama_thm} and defining
$$
G_n =\sum_{k=0}^{n+1}(-4)^k\frac{B_{2k} B_{2n+2-2k}}{(2k)!(2n+2-2k)!} ,
$$
% \begin{multline*}
% \frac{\pi^{-n}}{2^{-n}}\left\{\frac{1}{2}\zeta(2n+1)+S_{2n+1}(\pi) \right\} = \\ {(-2\pi)}^{-n}\left\{\frac{1}{2}\zeta(2n+1)+S_{2n+1}(4\pi) \right\} - 2^{-n-1}4^n\pi^{n+1}\sum_{k=0}^{n+1}(-4)^k\frac{B_{2k} B_{2n+2-2k}}{(2k)!(2n+2-2k)!}
% \end{multline*}
% 
% \begin{multline*}
% -2^n\left\{\frac{1}{2}\zeta(2n+1)+S_{2n+1}(\pi) \right\} = \\ 2^{-n}\left\{\frac{1}{2}\zeta(2n+1)+S_{2n+1}(4\pi) \right\} +2^{-n-1}4^n\pi^{2n+1}\sum_{k=0}^{n+1}(-4)^k\frac{B_{2k} B_{2n+2-2k}}{(2k)!(2n+2-2k)!}
% \end{multline*}
% 
% 
% Now we simplify:
% \begin{multline*}
% -2^n\left\{\frac{1}{2}\zeta(2n+1)+S_{2n+1}(\pi) \right\} = \\ 2^{-n}\left\{\frac{1}{2}\zeta(2n+1)+S_{2n+1}(4\pi) \right\} +2^{-n-1}4^n\pi^{2n+1}G_n
% \end{multline*}
% 
% \begin{equation*}
% -4^n\left\{\frac{1}{2}\zeta(2n+1)+S_{2n+1}(\pi) \right\} = \left\{\frac{1}{2}\zeta(2n+1)+S_{2n+1}(4\pi) \right\} + \frac{4^n}{2}\pi^{2n+1}G_n
% \end{equation*}
% 
% \begin{equation*}
% \frac{-4^n}{2}\zeta(2n+1)-4^nS_{2n+1}(\pi)  = \frac{1}{2}\zeta(2n+1)+S_{2n+1}(4\pi)  + \frac{4^n}{2}\pi^{2n+1}G_n
% \end{equation*}
% 
% \begin{equation*}
% \frac{-4^n}{2}\zeta(2n+1)  = \frac{1}{2}\zeta(2n+1)+S_{2n+1}(4\pi) +4^nS_{2n+1}(\pi) + \frac{4^n}{2}\pi^{2n+1}G_n
% \end{equation*}
% 
% 
% \begin{equation*}
% \frac{-(4^n+1)}{2}\zeta(2n+1)  = S_{2n+1}(4\pi) +4^nS_{2n+1}(\pi) + \frac{4^n}{2}\pi^{2n+1}G_n
% \end{equation*}
% 
one has
\[
\zeta(2n+1)  
= - \frac{2 \cdot 4^nS_{2n+1}(1)  + 2 S_{2n+1}(4) 
+  4^n\pi^{2n+1}G_n}{4^n+1} .
\]
Combining this with equation (\ref{eqn_(3)}) yields

% \begin{equation*}
% \frac{S_{2n+1}(4\pi) +4^nS_{2n+1}(\pi) + \frac{4^n}{2}\pi^{2n+1}G_n}{-(4^n+1)}+S_{2n+1}(2\pi)  = \frac{- 4^n}{2}\pi^{2n+1}F_n
% \end{equation*}
% 
% \begin{equation*}
% S_{2n+1}(4\pi) +4^nS_{2n+1}(\pi) + \frac{4^n}{2}\pi^{2n+1}G_n-(4^n+1)S_{2n+1}(2\pi)  = \frac{4^n}{2}\pi^{2n+1}(4^n+1)F_n
% \end{equation*}
% 
% \begin{equation*}
% 4^nS_{2n+1}(\pi) -(4^n+1)S_{2n+1}(2\pi) +S_{2n+1}(4\pi)= \frac{4^n}{2}\pi^{2n+1}(4^n+1)F_n-\frac{4^n}{2}\pi^{2n+1}G_n
% \end{equation*}
% 
% \begin{equation*}
% 4^nS_{2n+1}(\pi) -(4^n+1)S_{2n+1}(2\pi) +S_{2n+1}(4\pi)= \pi^{2n+1}\left(\frac{ 4^n}{2}(4^n+1)F_n-\frac{4^n}{2}G_n\right)
% \end{equation*}
% 

$$
\frac{4^nS_{2n+1}(1) -(4^n+1)S_{2n+1}(2) +S_{2n+1}(4)}{\frac{4^n}{2}(4^n+1)F_n-\frac{4^n}{2}G_n}= \pi^{2n+1}.
$$
Substituting $n=2m-1$  and defining
$$
D_m = \frac{4^{2m-1} \big[ (4^{2m-1}+1)F_{2m-1} - G_{2m-1} \big]}{2}
$$
produces

% \begin{equation*}
% \frac{4^{2m-1}S_{4m-1}(\pi) -(4^{2m-1}+1)S_{4m-1}(2\pi) +S_{4m-1}(4\pi)}{\frac{4^{2m-1}}{2}\{(4^{2m-1}+1)F_n-G_n\}}= \pi^{4m-1}
% \end{equation*}
% 

$$
\pi^{4m-1} = \frac{4^{2m-1}}{D_m} S_{4m-1}(1) -
\frac{4^{2m-1}+1}{D_m} S_{4m-1}(2) + \frac{1}{D_m} S_{4m-1}(4). 
$$
This identity may be combined with equation (\ref{eqn_(3)})
to give

$$
\zeta(4m-1) = 
-\frac{F_{2m-1} 4^{4m-2}  }{D_m} S_{4m-1}(1)
+\frac{G_{2m-1} 4^{2m-1}  }{D_m} S_{4m-1}(2)
-\frac{F_{2m-1} 4^{2m-1}  }{D_m} S_{4m-1}(4) .
$$

% \begin{multline*}
% -\frac{2F_n4^nS_{2n+1}(\pi)}{(4^n+1)F_n-G_n}+\frac{ 2{G_n}S_{2n+1}(2\pi) }{(4^n+1)F_n-G_n}
% -\frac{2{F_n}S_{2n+1}(4\pi)}{{(4^n+1)F_n-G_n}} =
% \\ \zeta(2n+1) 
% \end{multline*}
% 
% And the desired expression follows from substituting $n=2m-1$, so that $2n+1$ becomes $4m-1$.
% 

To obtain formulas for the $4m+1$ cases, set  $\alpha=2 \pi$,  $\beta=\pi/2$, and $n=2m$ in Theorem \ref{Rama_thm}
to obtain
\begin{equation}
\label{eqn_(4)}
\zeta(4m+1)  =  \frac{- 2 \cdot 16^m S_{4m+1}(1) + 2 S_{4m+1}(4)  
-  16^m \pi^{4m+1}G_{2m}}{16^m-1}.
\end{equation}
Define $T_n (r)$ by
$$
T_{n}(r)=\sum_{k=1}^{\infty}\frac{1}{k^{n}(e^{\pi r k}+1)},
$$
and another finite sum of Bernoulli numbers by
$$
H_{n}=\sum_{k=0}^{n}(-4)^{n+k}\frac{B_{4k} B_{4n+2-4k}}{(4k)!(4n+2-4k)!} .
$$

We begin with the case $m \ge 1$. Under this hypothesis, 
Vepstas  established the following expression (refer to the calculation following \cite[Theorem 7]{Vepstas}, and the statement in the introduction of \cite{Vepstas}):
\begin{eqnarray*}
\label{eqn_Vepstas}
\lefteqn{ (1+(-4)^m-2^{4m+1}) \zeta(4m+1)= } \\
& & 
2T_{4m+1}(2)+2(2^{4m+1}-(-4)^m)S_{4m+1}(2)
+2^{4m+1}\pi^{4m+1}H_m+2^{4m}\pi^{4m+1}G_{2m} .
\end{eqnarray*}
Vepstas also gave  a formula relating
$T_k$ and $S_k$:

$$
T_k(x)=S_k(x)-2S_k(2x).
$$
Combining the last two equations produces

\begin{eqnarray*}
\lefteqn{
\frac{1+(-4)^m-2^{4m+1}}
{\frac{1}{2}(1-4^{2m})}
\left(\frac{4^{2m}}{2}\pi^{4m+1}G_{2m}-S_{4m+1}(4)+4^{2m}S_{4m+1}(1)\right)
= }\\
& & 
2[2^{4m+1}-(-4)^m+1]S_{4m+1}(2)
-4S_{4m+1}(4)+2^{4m+1}\pi^{4m+1}H_m+2^{4m}\pi^{4m+1}G_{2m}.
\end{eqnarray*}
Letting
$$
K_m=\frac{\frac{1}{2}(1-4^{2m})}{1+(-4)^m-2^{4m+1}}
$$
and
$$
E_m = \frac{4^{2m}}{2}G_{2m}-2^{4m+1}K_mH_m-2^{4m}K_mG_{2m} ,
$$
one finds
$$
\pi^{4m+1} =
-\frac{4^{2m}}{E_m} S_{4m+1}(1)
+\frac{ 2K_m[2^{4m+1}-(-4)^m+1]}{E_m} S_{4m+1}(2)
+  \frac{ (1-4K_m)}{E_m}S_{4m+1}(4).
$$
Substituting this into equation (\ref{eqn_(4)}) produces

\begin{eqnarray*}
\zeta(4m+1) & = &
-\frac{ 16^m( 2E_m  - 16^m G_{2m} ) }{(16^m - 1)E_m} S_{4m+1}(1) \\
& &
-\frac{ 2\cdot 16^m G_{2m}K_m  (2\cdot 16^m - (-4)^m + 1) }
  {(16^m - 1)E_m} S_{4m+1}(2) \\ 
& & 
-\frac{16^m  G_{2m}  (1-  4 K_m ) -  2E_m }
  {(16^m - 1)E_m} S_{4m+1}(4) .
\end{eqnarray*}

It remains to consider the case $m=0$ and establish the formula for $\pi$ stated in the introduction. This formula is a consequence of classical results about $q$-series, which permit the direct evaluation of $S_1(1)$, $S_1(2)$, and $S_1(4)$. First, note that \cite[Equation (22.11)]{Berndt} gives 
\[
S_1(2) = \frac{1}{4} \log\left( \frac{4}{\pi} \right) - \frac{\pi}{12} + \log \Gamma \left( \frac{3}{4} \right),
\]
where $\log$ denotes the natural logarithm. A straightforward adaptation of the proof of this formula  in \cite{Berndt} yields, after using the tables of Zucker cited therein with the choice $c^2 = 4$,
\[
S_1(4) =  - \frac{1}{6} \log ( 2^{-33/4} ) + \frac{1}{6} \log ( \pi^{9/2} ) - \log \Gamma \left( \frac{1}{4} \right)  - \frac{\pi}{6}.
\]
Further, \cite[Equation (2.4)]{berndt2017ramanujan} gives
\[
S_1 (1) = S_1(4) + \frac{1}{4} \log \left( \frac{1}{4} \right)  + \frac{\pi}{8}.
\]
Substituting these expressions into the claimed formula for $\pi$ and simplifying completes the proof.

%%%%%%%%%%%%%%%%%%%%%%%%%%%%%%%%%%%%%%%%%%%%%%%%%%%%%%%%%%%%%%%%%

\bigskip
\hrule
\bigskip

\noindent 2010 {\it Mathematics Subject Classification}:
Primary 11Y60.

\noindent \emph{Keywords: }
$\pi$, Riemann zeta function.

\bigskip
\hrule
\bigskip

\vspace*{+.1in}
\noindent
Received January 24 2011;
revised version received  February 7 2011.
Published in {\it Journal of Integer Sequences}, February 20 2011.
ArXiv version (revised and corrected) posted June 1, 2024. 
\bigskip
\hrule
\bigskip

%\noindent
%Return to
%\href{Journal of Integer Sequences home page}{\tt  http://www.cs.uwaterloo.ca/journals/JIS/}.
%\vskip .1in

\end{document}